\documentclass{amsart}

\usepackage{graphicx}%
\usepackage{multirow}%
\usepackage{amsmath,amssymb,amsfonts}%
\usepackage{amsthm}%
\usepackage{mathrsfs}%
\usepackage[title]{appendix}%
\usepackage{xcolor}%
\usepackage{textcomp}%
\usepackage{manyfoot}%
\usepackage{booktabs}%
\usepackage{algorithm}%
\usepackage{algorithmicx}%
\usepackage{algpseudocode}%
\usepackage{listings}%

\newtheorem{theoremv}{Main Theorem}

\theoremstyle{thmstyledefinition}%
\newtheorem*{defi}{Definition}

\newcommand{\block}[2]{\addtolength{\arraycolsep}{-3pt}\begin{pmatrix}#1\cr#2\end{pmatrix}}

\def\C{\mathbb{C}}
\def\cat{\mathop{\mathrm{cat}}}
\def\H{\mathbb{H}}
\def\I{\mathbf{i}}
\def\ii{\I}
\def\J{\mathbf{j}}
\def\K{\mathbf{k}}
\def\Log{{\mathrm{Log}\,}}
\def\mega{\mathbf{\omega}}
\def\Pv{{\mathrm{P}}}
\def\R{\mathbb{R}}
\def\Sp{{\mathrm{Sp}}}
\def\su{{\mathrm{SU}}}

\begin{document}

\title[Subspace LS-category of quasi-projective quaternionic spaces]{Subspace {L}usternik-{S}chnirelmann category of quasi-projective quaternionic spaces}

\author[E. Mac\'ias-Virg\'os --- D. Tanr\'e]{Enrique Mac\'ias-Virg\'os --- Daniel Tanr\'e}



\thanks{The two authors are  supported  by the MICINN research project 
PID2020-114474GB-100 and the ANR-11-LABEX-0007-01  ``CEMPI''.
The first author is partially supported by 
Consolidaci\'on 2023 ED431C 2023/31 GiMAT, Proxectos Plan galego de financiamento universitario para o periodo 2022-2026.}

\begin{abstract}%
Let $Q_n$ be the quasi-projective subspace of the symplectic group $\mathrm{Sp}(n)$.
In this short note, we prove that the subspace Lusternik-Schnirelmann category of $Q_n$ in $\mathrm{Sp}(n)$ is 2.
For that, we use a quaternionic logarithm, as Singhof did in the complex case for the 
determination of the Lusternik-Schnirelmann category of the unitary group.
Our result generalizes the known case $n=2$ (by L. Fern\'andez-Su\'arez, A. G\'omez-Tato and D. Tanr\'e) and has to be compared to the equality $\mathrm{cat}\,Q_{3}=3$, established by N. Iwase and T. Miyauchi.
\end{abstract}

\maketitle



\section*{Introduction}

The problem list of T. Ganea (\cite{MR0339147}) still contains unresolved questions on the determination of 
Lusternik-Schnirelmann category (in short, LS-category).
In particular, the LS-category of classical Lie groups is still incomplete.
The complex case was completely solved by W.~Singhof (\cite{Singhof1975a,Singhof1976}) with $\cat\su(n)=n-1$.
In the quaternionic case, the LS-category of the \emph{symplectic group}, $\Sp(n)$, is known only in low dimensions: 
$\cat\,\Sp(1)=\cat \,S^3=1$,
$\cat\,\Sp(2)=3$ (\cite{Schw1965})
and
$\cat\,\Sp(3)=5$ (\cite{FerGomTanStr2004}).
For the higher cases, there are few informations. One knows from \cite{IwaMim01} that $\cat\,\Sp(n)\geq n+2$ if $n\geq 3$,
and from \cite{MacPer2013} that $\cat\,\Sp(n)\leq n(n+1)/2$.

In this work, we are interested in quasi-projective spaces. 
In the complex case, they coincide with the suspension
of the projective spaces $\C\Pv^n$ but, in the quaternionic case, the situation is different. 
Their cellular structure still determines the cohomology ring of the group $\Sp(n)$ (\cite{Steenrod1962})
but they are not the third suspension of the quaternionic projective spaces $\H\Pv^n$.
Denote them by $Q_{n}$ and look at the first cases.

For $n=1$, we have the sphere $S^3$. For $n=2$ we get a 2-cell space $Q_{2}=S^3\cup_{\varphi}e^7$ where 
$\varphi$ is a generator of $\pi_{6}(S^3)$. It has a non-zero Hopf invariant, thus $\cat Q_{2}=2$,
see \cite[Chapter 6]{book}.
The first difficult case, $n=3$, was solved by N. Iwase and  T. Miyauchi (\cite{IwaMiy2007}) who proved
$\mathrm{cat}\,Q_{3}=3$.

In contrast with this result, let us notice that  \cite[Theorem 3.3]{FerGomTan2001} implies 
that the subspace $Q_{3}$ has subspace category 2 in $\Sp(3)$.
The proof relies on the knowledge of some homotopy groups of $Q_{2}$, $\Sp(2)$, $\Sp(3)$.
We prove here that this is a general fact and the proof can be done 
``\`a la Singhof''
from the  description
of $Q_{n}$ and $\Sp(n)$ as sets of particular matrices and the use of a quaternionic logarithm.

\begin{theoremv}\label{prop:subcat2}
The subspace LS-category of $Q_{n}$ in $\Sp(n)$ is 2.
\end{theoremv}

Let us observe that  the product cells of $\su(n)$ are grafted onto the quasi-projective complex space
of LS-category 1. In the case of $\Sp(n)$, the quasi-projective quaternionic space is not a suspension.
Moreover it is not even a space of constant LS-category for any $n$, as shown by the cases $n=2$ and $n=3$.
But, for a determination of the LS-category of $\Sp(n)$, 
the main point is not  the LS-category of the quasi-projective space
but its relative LS-category in $\Sp(n)$. 
Our result shows that the two situations are similar.
It reinforces the intuition of the existence of an upper bound for the LS-category  of $\Sp(n)$ 
in the form of a polynomial in $n$ of degree 1, the most challenging  conjecture being 
$\cat \Sp(n)\leq 2n-1$.

\medskip

The basic definitions---such as those of quasi-projective spaces, 
subspace LS-category, quaternionic logarithm, and so on---are recalled in the first section. 
The second section  contains the proof of the Main theorem.

\section{Background}\label{sec:back}

Let us begin with the \emph{symplectic group}, $\Sp(n)$. We denote by
$\H$  the (non-commutative) field of quaternions and  $\H^n$ the $n$-dimensional vector space over $\H$,
with a right $\H$-action.
If $x=(x_{i})_{1\leq i\leq n}\in \H^n$ and $y=(y_{i})_{1\leq i\leq n}\in \H^n$, we define the 
scalar product 
$\langle x,y\rangle=\sum_{i=1}^n\overline{x}_{i}y_{i}$,
where $\overline{x}_{i}$ is the conjugate of $x_{i}\in\H$.
The  group $\Sp(n)$ is the group of orthogonal transformations of $\H^n$ or, equivalently, of the  quaternionic matrices  of order $n$ such that
$A^*A=I_{n}$.
We embedd $\H^n$ in $\H^{n+1}$ by putting the last coordinate equal to zero,
and $\Sp(n)$ in $\Sp(n+1)$ by
$$A\mapsto
\begin{pmatrix}
A&0\cr
0&1\cr
\end{pmatrix}.
$$

The use of the logarithm in the study of LS-category of Lie groups takes its origin in the determination
 $\cat\su(n)=n-1$ made by Singhof in \cite{Singhof1975a,Singhof1976}. Let us recall the definition of a principal branch of the logarithm for quaternions.

\medskip
The Lie group of unit quaternions $S^3$ has for Lie algebra  the set $\H_0=\langle\I,\J,\K\rangle$ 
of skew-symmetric quaternions $q$; i.e., those of the form $t\,\mega$ for some imaginary unit quaternion 
$\mega$ 
and some real number $t\in\R$. 
Any quaternion $q\in \H$  can be written 
%
$$q=s+t \,\mega, \quad s,t\in\R,  t\geq 0, \quad s^2+t^2=\vert q \vert^2.$$
By taking the argument $\theta\in  [0,+\pi)$, we get the polar form
$$q= \vert q \vert e^{\theta\,\mega}=\vert q \vert\left(\cos\theta + (\sin\theta)\, \mega\right).$$
The \emph{principal branch of the quaternionic logarithm} of $q=s+t\mega$ is then defined  (\cite{GPV}) as
$$\Log q = \ln \vert q \vert +   \theta\,\mega=\ln\vert q \vert + \arccos \Re\left(\frac{q}{\vert q \vert}\right)\,\mega,$$
where
$\theta\in [0,+\pi)$ and $\Re$ is the real part.
This function is continuous except for 
$\theta=\pi$, which corresponds to the quaternions
$q=-\vert q \vert $
on the negative real axis. 
Let us notice that two similar quaternions  (i.e. same norm, same real part) have similar logarithms.

\medskip
Let $S^{4n-1}$ be the unit sphere of $\H^n$. The quaternionic quasi-projective space, $Q_{n}$, is the image of the map
$$\varphi\colon S^{4n-1}\times S^3\to \Sp(n),$$
defined  by
$$\varphi(x,\lambda)=x(\lambda -1)x^*+I_{n}.$$
(More geometrically, the map $\varphi(x,\lambda)$ sends $x$ to $x\lambda$
and is the identity on the  subspace orthogonal to $x$.)
Equivalently,  $Q_{n}$ is the quotient space of $S^{4n-1}\times S^3$ by the equivalence relation
$$(x,\lambda)\sim (x\nu,\nu^{-1}\lambda\nu), \text{ with } \nu\in S^3, \text{ and }
(x,1)\sim (y,1)\text{ for any } x,\,y\in S^{4n-1}.$$
%
%
The group law $\mu$ of $\Sp(n)$ induces a relative homeomorphism (\cite[Proposition 2.4]{Steenrod1962})
$$\mu\colon (Q_{n}\times \Sp(n-1),Q_{n-1}\times \Sp(n-1))
\to
(\Sp(n),\Sp(n-1)),$$
which maps $Q_{n}\times \Sp(n-1)$ onto $\Sp(n)$.

Let $E^{4(n-1)}$ be the ball consisting of all vectors $x\in S^{4n-1}\subset \H^n$ with $x_{n}\in \R$
and $x_{n}\geq 0$. We denote  $h_{n}$ the composition of $\varphi$ with the product of the inclusion map
$\iota_{n}\colon E^{4(n-1)}\to S^{4n-1}$
and the canonical surjection $\rho\colon E^3\to S^3\cong E^3/S^2$,
$$h_{n}\colon E^{4n-1}=E^{4(n-1)}\times E^3\xrightarrow{\iota_{n}\times \rho}
S^{4n-1}\times S^3 \xrightarrow{\varphi}Q_{n}.$$
The characteristic maps of the cells of $\Sp(n)$ are products of maps $h_{i}$,
$$E^{4i_{1}-1}\times\cdots\times E^{4i_{r}-1}
\xrightarrow{h_{i_{1}}\times\cdots\times h_{i_{r}}}
Q_{i_{1}}\times\cdots\times Q_{i_{r}}
\xrightarrow{\mu}
\Sp(n),$$
where $n\geq i_{1}>i_{2}>\dots >i_{r}>0$. Such a cell is denoted by $(i_{1},\dots,i_{r})$ and called
a normal cell. In \cite[Theorem 2.1]{Steenrod1962}, Steenrod shows that  
 $\Sp(n)$ is a CW-complex whose cells are the normal cells and the 0-cell $I_{n}$.
 \medskip
\begin{defi}
 The \emph{subspace category}
 of a subspace $A$ of a topological space $X$, denoted $\cat_{X}A$, is the least integer $n$ such that 
there exist open subsets $(U_{i})_{0\leq i\leq n}$ of $X$ which cover $A$ and which are contractible \emph{in} $X$.
\end{defi}
If $X$ is a metric ANR, one gets the same integer if the $U_{i}$'s are any subsets of $X$, 
see \cite[Theorem 3.2]{MR3193429}.
It is true for instance for any CW-complex.

\section{Proof of the Main Theorem}\label{sec:qn=2}

We need three subsets whose union is $Q_{n}$, each one being contractible in $\Sp(n)$. We set:
\begin{itemize}
\item $O_{1}=\left\{\varphi(x,\lambda)\mid x\in S^{4n-1},\,\lambda\neq 1,\,\lambda\neq -1\right\}$,
\item $O_{2}=\left\{\varphi(x,-1)\mid  x\in S^{4n-1}\right\}$ and
\item $O_{3}$ an open ball containing $I_{n}$.
\end{itemize}
By definition, they cover $Q_{n}$ and the open set $O_{3}$ is contractible. We now use the quaternionic logarithm to obtain an explicit formula for the contracting homotopies
in $O_{1}$ and $O_{2}$.

Let us begin with 
$A=\varphi(x,\lambda)=x(\lambda -1)x^*+I_{n}\in O_{1}$.
We choose an orthonormal basis $(x,Y)$ of $\H^n$, with $Y=(y_{2},\dots,y_{n})\in \H^{n-1}$, such that
$$A=\left(x \;Y\right)
\begin{pmatrix}\lambda&0\\
0&I_{n-1}
\end{pmatrix}
\block{x^*}{Y^*}.
$$
As $\lambda\neq -1$ the quaternionic logarithm is defined by 
$\Log A  = x (\Log \lambda) x^*$
and we have a continuous map, from $O_{1}$ to $\Sp(n)$,
$$\Log A = \left(x \;Y\right)
\begin{pmatrix}
\Log\lambda&0\\
0&I_{n-1}
\end{pmatrix}
\block{x^*}{Y^*}.
$$
By using  the exponential,  we get
\begin{eqnarray*}
A=e^{\Log A}
& = &
e^{x(\Log\lambda)x^*} 
\\
&=&
1+\sum_{n=1}^{\infty} \frac{1}{n!} (x(\Log \lambda)x^*)^n =
1+x\left(\sum_{n=1}^{\infty} \frac{1}{n!} (\Log \lambda)^n\right)x^*\\
&=&
1+x\left(e^{\Log\lambda}-1\right)x^*=1+x(\lambda -1)x^*=\varphi(x,\lambda).
\end{eqnarray*}
We define a continuous map from $O_{1}\times [0,1]$ to $\Sp(n)$ by
$$\Phi(A,t)=e^{x ((1-t)\Log\lambda) x^*}.$$
The extremities of this path are $\Phi(A,0)=e^{\Log A}=A$ and $\Phi(A,1)=I_{n}$. 
A computation similar to the previous one also gives
$$\Phi(A,t)=1+x\left(e^{(1-t)\Log\lambda}-1\right)x^*=\varphi(x,e^{(1-t)\Log\lambda})$$
and the contractibility of $O_{1}$.

\medskip
We continue with  $A=\varphi(x,-1)\in O_{2}$.  From the definition, we have
$\varphi(x,-1)=\varphi(x,\ii)\varphi(x,\ii)$. We define a continuous map
$\Psi\colon O_{2}\times [0,1]\to \Sp(n)$ by
$$\Psi(A,t)=\left\{
\begin{array}{ccl}
\varphi(x,e^{(1-2t)\Log\ii})\varphi(x,\ii)&\quad\text{if}&0\leq t\leq 1/2,\\
\varphi(x,e^{2(1-t)\Log \ii})&\quad\text{if}&1/2\leq t\leq 1.
\end{array}\right.$$
We observe that
$\Psi(A,0)=\varphi(x,\ii)\varphi(x,\ii)=\varphi(x,-1)$
and 
$\Psi(A,1)=\varphi(x,1)=I_{n}$.

\medskip
In the previous proof, the subsets  $O_{1}$ and $O_{3}$ are contractible and not only contractible in $\Sp(n)$. 
In contrast, this is not the case for $O_{2}$: this open subset is contractible in $\Sp(n)$ but the contraction
does not stay in $Q_{n}$. This is coherent with the fact that the LS-category of $Q_{n}$ is greater than 3 
(\cite{IwaMiy2007}) for $n\geq 3$.

\bigskip

{\textsc{Enrique Mac\'ias-Virg\'os}\\
CITMAga, Departamento de Matem\'aticas, \\
Universidade de Santiago de Compostela, \\
15782 Santiago de Compostela, Spain}\\
\email{quique.macias@usc.es}\\

{\textsc{Daniel Tanr\'e}\\
D\'epartement de Math{\'e}matiques,\\ UMR-CNRS 8524, Universit\'e de Lille,\\
 59655 Villeneuve d'Ascq Cedex, France}\\
\email{Daniel.Tanre@univ-lille.fr}


\begin{thebibliography}{99}
\labelsep=1em\relax

\bibitem{book}
    O. Cornea, 
    G. Lupton, 
    J. Oprea, 
    D. Tanré:  
    Lusternik-Schnirelmann Category.  
    Mathematical Surveys and Monographs, vol. 103.  
    American Mathematical Society, Providence, RI (2003)

\bibitem{FerGomTanStr2004}
    L. Fernández-Suárez, 
    A. Gómez-Tato, 
    J. Strom, 
    D. Tanré:  
    The Lusternik-Schnirelmann category of $\Sp(3)$.  
    Proc. Amer. Math. Soc. 132(2), 587--595 (2004)

\bibitem{FerGomTan2001}
    L. Fernández-Suárez, 
    A. Gómez-Tato, 
    D. Tanré:  
    Hopf-Ganea invariants and weak LS category.  
    Topology Appl. 115(3), 305--316 (2001)

\bibitem{MR0339147}
    T. Ganea:  
    Some problems on numerical homotopy invariants.  
    In: Symposium on Algebraic Topology (Battelle Seattle Res. Center, Seattle, Wash., 1971).  
    Lecture Notes in Math., vol. 249, pp. 23--30.  
    Springer, Berlin (1971)

\bibitem{GPV}
    G. Gentili, 
    J. Prezelj, 
    F. Vlacci:  
    On a definition of logarithm of quaternionic functions.  
    J. Noncommut. Geom. 17(3), 1099--1128 (2023)

\bibitem{IwaMim01}
    N. Iwase, 
    M. Mimura:  
    L-S categories of simply-connected compact simple Lie groups of low rank.  
    In: Categorical Decomposition Techniques in Algebraic Topology (Isle of Skye, 2001).  
    Progr. Math., vol. 215, pp. 199--212.  
    Birkhäuser, Basel (2004)

\bibitem{IwaMiy2007}
    N. Iwase, 
    T. Miyauchi:  
    Lusternik-Schnirelmann category of stunted quasi-projective spaces.  
    J. Math. Kyoto Univ. 47(2), 321--326 (2007)

\bibitem{MacPer2013}
    E. Mac\'ias-Virg\'os,
    M.J. Pereira-S\'aez,
    An upper bound for the Lusternik-Schnirelmann category of the symplectic group.  
    Math. Proc. Cambridge Philos. Soc. 155(2), 271--276 (2013)

\bibitem{Schw1965}
	P.A. Schweitzer:  
    Secondary cohomology operations induced by the diagonal mapping.  
    Topology 3, 337--355 (1965)

\bibitem{Singhof1975a}
    W. Singhof:  
    On the Lusternik-Schnirelmann category of Lie groups.  
    Math. Z. 145(2), 111--116 (1975)

\bibitem{Singhof1976}
    W. Singhof:  
    On the Lusternik-Schnirelmann category of Lie groups. II.  
    Math. Z. 151(2), 143--148 (1976)

\bibitem{MR3193429}
    T. Srinivasan:  
    The Lusternik-Schnirelmann category of metric spaces.  
    Topology Appl. 167, 87--95 (2014)

\bibitem{Steenrod1962}
	N.E. Steenrod:  
    Cohomology Operations.  
    Annals of Mathematics Studies, No. 50.  
    Princeton University Press, Princeton, N.J. (1962).  
    Lectures by N.E. Steenrod written and revised by D.B.A. Epstein.

\end{thebibliography}
\end{document}